\documentclass[12pt]{amsart}
\usepackage{amsthm,amsmath,amssymb,amsfonts,amscd,verbatim,latexsym}
\newtheorem{Theorem}{Theorem}[section]
\newtheorem{Lemma}[Theorem]{Lemma}
\newtheorem{Proposition}[Theorem]{Proposition}

\newtheorem{Remark}[Theorem]{Remark}

\newcommand{\demo}{\noindent {\sc Proof.\;}}
\renewcommand{\P}{{\mathbb P}} 
\renewcommand{\O}{{\mathcal O}} 
\newcommand{\E}{{\mathcal E}} 
\renewcommand{\L}{{\mathcal L}} 
\newcommand{\vspan}{\text{span}}
\newcommand{\red}{\text{red}}
\newcommand{\Hilb}{\text{Hilb}}
\newcommand{\ra}{\rightarrow}

\newcommand{\lra}{\longrightarrow}

\newcommand{\II}{\mathcal I}

\newcommand{\rank}{\text{rank}\,}
\newcommand{\Proj}{\text{Proj}\,}

\newcommand{\complex}{{\mathbf C}}

\newcommand{\codim}{\text{codim}\,}

\newcommand{\U}[1]{{\mathcal U}_{\, #1}}

\newcommand{\brac}[3]{(#1 \, #2 \, #3)}
\begin{document}
\title{The Waring loci of ternary quartics} 
\author{Jaydeep V.~Chipalkatti} 
\maketitle 

\parbox{12cm}{\small 
{\sc Abstract.} 
Let $s$ be any integer between $1$ and $5$. We determine necessary 
and sufficient conditions that a given ternary quartic be expressible 
as a (possibly degenerate) sum of fourth powers of $s$ linear 
forms.}

\vspace{2mm}

\mbox{\small AMS subject classification: 14-04, 14 L35.}

\mbox{\small Keywords: apolarity, Waring's problem, concomitants.}

\bigskip \bigskip 

Let $F$ be a homogeneous polynomial (or form) 
of degree $q$ in $r$ variables. 
It is a classical problem to decide whether $F$ can be expressed as a 
sum of powers of linear forms, 
\begin{equation} F = L_1^q + \dots + L_s^q ; 
\label{eq.Waring} \end{equation}
for a specified number $s$. This is usually called `Waring's problem' 
for algebraic forms. 
(Normally one also allows a `degeneration' of the right hand side in 
(\ref{eq.Waring}), this will be made precise later.) Since the condition 
that $F$ be so expressible is invariant under the natural action of 
the group $SL_r$, it should be equivalent to the vanishing of 
certain concomitants of $F$. It is of interest to identify these 
concomitants, and thus to get explicit algebraic conditions on 
the coefficients of $F$ for the expression (\ref{eq.Waring}) (or 
its degeneration) to be possible. 

In this note we consider ternary quartics, i.e., we let $q=4, r=3$. 
Since a general ternary quartic is a sum of $6$ powers of linear 
forms, we consider the range $1 \le s \le 5$. 
The calculations required in this case are not prohibitively 
large, and it is possible to get a complete solution. 
The result is given in Theorem \ref{main.theorem}. 

For an excellent introduction to Waring's problem, see \cite{Ger1}. 
A very comprehensive account of the theory is given in 
\cite{IK}. The problem was solved for binary forms by 
Gundelfinger (see \cite{ego4, GrYo, Kung2}). For ternary cubics, the 
solution is essentially given in \cite{Salmon2}. 

\section{Preliminaries}
In this section, we establish notation and recall the 
representation-theoretic notions that we will need. The reader may 
wish to read it along with \cite{ego5}, where a similar 
set-up is used but more detailed explanations are given. 
Although we work throughout with ternary quartics, on occasion
I state whether a result goes through for arbitrary $q$ and $r$. 

Let $V$ be a three-dimensional $\complex$-vector space. 
Let $V$ and $V^*$ have dual bases $\{y_0,y_1,y_2\}$ and 
$\{x_0, x_1, x_2 \}$ respectively. We will identify 
$\P \, ( S_4 \, V^*) = \P^{14}$ with the space of 
quartic forms in the $x_i$ (up to scalars). Let 
$R = S_\bullet (S_4 \, V)$ denote the symmetric algebra on 
$S_4 \, V$, then $\P^{14} = \text{Proj} \, R$. 

For an integer $s$, let 
$W_s^\circ \subseteq \P^{14}$ denote the set 
\begin{equation}
\{ F \in \P^{14}: 
F = L_1^4 + \dots + L_s^4,  \; \text{for some $L_i \in V^*$} \} , 
\label{eq.Ws} \end{equation}
and $W_s$ its Zariski closure. The $W_s$ are irreducible projective 
varieties, whereas the $W_s^\circ$ are generally 
not quasiprojective, but only constructible. The 
imbedding $W_s \subseteq \P^{14}$ is $SL(V)$-equivariant. 

We would like to give necessary and sufficient $SL_3$-invariant 
algebraic conditions for a form $F$ to lie in $W_s$. This is 
roughly the same as determining the structure of the ideal 
$I_{W_s} \subseteq R$ as an $SL(V)$-representation. The parallel 
problem of characterizing $W_s^\circ$ is harder, and later in 
\S \ref{w2open.subsection} we only consider the case 
of $W_2^\circ$. 

\begin{Remark} \rm 
The dimensions of the $W_s$ are known by a theorem due to 
Alexander and Hirschowitz - see \cite{Ger1,IK} for references. In fact 
$W_s = \P^{14}$ for $s \ge 6$, $\dim W_s = 3s-1$ for $1 \le s \le 4$ 
and $W_5$ is a hypersurface. The degrees of $W_1, \dots, W_5$ are 
$16, 75, 112, 35$ and $6$ respectively (see \cite{Norway1}). 
Elements of $W_4$ (resp.~$W_5$) are conventionally 
called Capolary (resp.~Clebsch) quartics. 
\end{Remark} 

\subsection{Schur functors} 
If $\lambda$ is a partition, then $S_\lambda(-)$ denotes the 
corresponding Schur functor. (We maintain the indexing conventions 
of \cite[Ch.~6]{FH}.) 

Let $J \subseteq R$ be a homogeneous $SL(V)$-stable 
ideal, and $J_d$ its degree $d$ part. We have an additive decomposition 
\begin{equation}
J_d = \bigoplus\limits_\lambda \, (S_\lambda V)^{N_\lambda} 
\subseteq S_d \, (S_4 \, V). 
\label{eq.schur} \end{equation}
Since $V$ is three-dimensional, $\lambda = (m+n,n)$ for some 
integers $m,n$. If we can locate the degrees $d$ which generate $J$, and 
then specify the inclusions in (\ref{eq.schur}), then $J$ is 
completely specified. These inclusions are encoded by the concomitants 
of ternary quartics. 

\subsection{Concomitants} 
Write $a_I$ for the monomial $y_0^{i_0}y_1^{i_1}y_2^{i_2}$, where 
$I = (i_0, i_1, i_2)$ is of total degree $4$. Then 
the $\{a_I\}$ form a basis of $S_4 \, V$, and $R$ is the polynomial 
algebra $\complex \, [\{a_I\}]$. Now 
$S_4 \, V \otimes S_4 \, V^*$ contains the trace element, which, 
when written out in full, appears as 
\[ {\mathbb F} = 
\sum\limits_{|I| = 4} a_I \otimes x_0^{i_0}x_1^{i_1}x_2^{i_2}.\] 
This is the precise formulation of the concept of a `generic' 
ternary quartic. Write 
\[ u_0 = x_1 \wedge x_2, \; u_1 = x_2 \wedge x_0, \; u_2 = x_0 \wedge x_1; 
\] 
which form a basis of $\wedge^2 V^*$. 
Consider an inclusion 
\[ S_{(m+n,n)} V \stackrel{\varphi}{\lra} S_d \, ( S_4 \, V) \] 
of $SL(V)$-representations, which will correspond to an inclusion 
\[ \complex \lra S_d \, (S_4 \, V) \otimes S_{m+n,n} V^*. \] 
Let $\Phi$ denote the image of $1$ under this map. Written out in 
full, it appears as a form 
of degree $d,m,n$ respectively, in the three sets of variables 
$a_I, x_i,u_i$. 
(This follows because $S_{m+n,n} \, V^*$ has a basis coming from 
standard tableaux, see \cite[Ch,~6]{FH}.) 
We will write this image as $\Phi(d,m,n)$. Classically, 
it is called a concomitant of degree $d$, order $m$ and class $n$ 
(of ternary quartics). For instance, ${\mathbb F}$ itself is 
a $\Phi(1,4,0)$ and the Hessian of $F$ is a $\Phi(3,6,0)$. 
For fixed $d,m,n$, the number of linearly independent 
concomitants equals the multiplicity of $S_{m+n,n}\, V$ in 
$S_d(S_4 \, V)$. 

We may regard $\Phi$ as a form in $u_i,x_i$ with coefficients 
in $R_d$. Then the subspace of $R_d$ generated by these coefficients 
coincides with the image of the inclusion $\varphi$ above. 
We can evaluate $\Phi$ at a specific form $F$ by substituting its 
actual coefficients for the letters $a_I$. Then we say that 
$\Phi$ vanishes at $F$ if this evaluated form is identically zero. 

\section{Apolarity}
We briefly explain the connection between the expression of 
$F$ as a sum of powers, and the existence of schemes `apolar' to 
$F$ (also see \cite{DolgachevKanev, EhRo,IK}). 

Let $A = S_\bullet V$ denote the symmetric algebra on $V$. 
A linear form $L \in V^*$ can be considered as a point in 
$\P V^* = \Proj A$. If $Z = \{L_1, \dots, L_s \}$ is a collection of 
points in $\P V^*$, then $I_Z \subseteq A$ denotes their ideal. 

For every $k \ge 0$, there is a coproduct map 
\[ S_4 \, V^* \lra S_k \, V^* \otimes S_{4-k} \, V^*. \] 
(It is zero unless $0 \le k \le 4$.) 
This gives rise to a map 
\begin{equation}
 S_4 \, V^* \lra \text{Hom}\, (S_k \, V, S_{4-k} \, V^*), \quad 
F \lra \alpha_{k,F}. 
\label{eq.hom} \end{equation}
Taking a direct sum over all $k$, we have a map 
\[  S_4 \, V^* \lra 
\text{Hom}\, (A , \, \bigoplus\limits_{i=0}^4 \, S_i \, V^*), \quad 
F \lra \alpha_F. \] 

For a fixed $F$, $\ker \alpha_F$ is 
a homogeneous ideal in $A$. Classically, a form in 
$\ker \alpha_F$ is said to be apolar to $F$. The passage between 
apolarity and Waring's problem is forged by the following beautiful 
theorem of Reye. (See \cite[Theorem 5.3 B]{IK}.)

\begin{Theorem}[Reye] \sl 
With notation as above, 
\[ 
F \in \vspan \, \{L_1^4, \dots, L_s^4 \} \iff 
I_Z \subseteq \ker \alpha_F \iff (I_Z)_4 \subseteq \ker \alpha_{4,F}. 
\] 
\end{Theorem}
The result holds for any $q, r$. In the situation above, we will 
say that $Z$ is an apolar scheme of $F$. 
To restate the theorem, $F$ lies in $W_s^\circ$ iff $F$ has a 
reduced zero-dimensional apolar scheme $Z$ of length $s$. 
It is natural to relax the requirement that $Z$ be reduced, which 
motivates the following definition. Let 
\[ X_s = \{ F \in \P^{14}: \ker \alpha_F \supseteq I_Z 
\; \; \text{for some $Z \in \text{Hilb}^s(\P V^*)$} \}.
\] 
This is a projective subvariety of $\P^{14}$, which a 
priori only contains $W_s$. Now \cite[Prop. 6.7]{IK} implies the 
following: 
\begin{Lemma} We have $W_s = X_s$ for all $s$. 
\end{Lemma} 

The result is true for $r = 2, 3$ and all $q$. 
(The crux is that the schemes 
$\Hilb^s(\P^{r-1})$ are irreducible for $r=2,3$.) I do not know 
if $X_s$ can properly contain $W_s$ for $r>3$. 

One can use Reye's theorem to relate the varieties 
$X_s$ to degeneracy loci of certain morphisms of vector bundles 
on $\P^{14}$. 
Globally, the map (\ref{eq.hom}) gives a morphism of vector bundles 
\[ \alpha_k: S_k \, V \otimes \O_{\P^{14}}(-1) \lra  S_{4-k} \,V^*.  \] 
Upto a twist, $\alpha_k$ is dual to $\alpha_{4-k}$. 

\begin{Lemma} \sl 
If $F$ belongs to $X_s$, then rank $\alpha_{k,F} \le s$ for any $k$. 
\end{Lemma}
\demo
Since $(I_Z)_k \subseteq \ker \alpha_{k,F}$, we have 
\[ \rank \alpha_{k,F} = \codim ( \ker \alpha_{k,F}, A_k ) 
\le \codim ( (I_Z)_k, A_k ) \le s. \] 

If $\psi: {\mathcal F} \lra \E$ is a morphism of vector bundles, 
let $Y(s, \psi)$ denote the scheme $\{ \rank \psi \le s\}$, 
whose ideal sheaf is locally generated by 
the $(s+1) \times (s+1)$-minors of a matrix representing $\psi$. 
Let $Y_\red(s, \psi)$ denote the underlying variety. 
We will shorten this to $Y$ or $Y_\red$ if no confusion is likely. 

\begin{Remark} \rm 
By the lemma above, $X_s \subseteq Y_\red(s, \alpha_k)$ 
for any $k$. For binary forms, the containment 
\[ X_s \subseteq \; \bigcap\limits_k Y_\red(s, \alpha_k) \] 
is an equality, but it is not known in general when this 
holds. 
\end{Remark} 

\subsection{Symmetric bundle maps} 
\label{symm.maps}
In the sequel, we will exploit the fact that $\alpha_{2,F}$ is a 
twisted symmetric morphism. 

Generally, let $T$ be a smooth complex projective 
variety. Let 
$\E$ be a rank $e$ vector bundle and $\L$ a line bundle 
on $T$. Assume that 
$ \psi: \E \lra \E^* \otimes \L $ is a twisted symmetric bundle map (i.e., 
$\psi^* \otimes \L = \psi$). 
Define $Y = Y(s,\psi)$ as above. Then assuming $Y$ is nonempty, 
\begin{equation} 
\codim (Y', T) \le \frac{(e-s)(e-s+1)}{2}, 
\label{eq.codim} \end{equation}
for every component $Y'$ of $Y$. 
Moreover, if equality holds for every component, then $Y$ is 
Cohen-Macaulay. In that case, the class of $Y$ in the Chow ring of $T$ 
is given by a determinantal formula (see \cite{HarrisTu}). 
Let $z_k = c_k(\E^* \otimes \sqrt{\L})$, then 
$[Y]$ equals $2^{e-s}$ times the $(e-s) \times (e-s)$ determinant whose 
$(i,j)$-th entry is $z_{(e-s-2i+j+1)}$. 

The minimal resolution of $Y$ (assuming equality in 
(\ref{eq.codim})) is deduced in \cite{JPW}. All that 
we need is the beginning portion 
\begin{equation}  
S_{\lambda_s} \E \otimes \L^{\otimes (-s-1)} \ra \O_T \ra \O_Y \ra 0, 
\label{res.jpw} \end{equation}
where $\lambda_s$ denotes the partition 
$(\underbrace{2, \dots, 2}_{s+1})$. 

We will apply the set-up to $\alpha_2$, with 
$\E = S_2 \, V \otimes \O_{\P^{14}}(-1)$ and 
$\L = \O_{\P^{14}}(-1)$.
\section{Computations}
We proceed to state the main theorem, and then explain the 
calculations entering into it. The concomitants will be written in 
the symbolic notation - see \cite{ego5,GrYo} for an 
explanation of this formalism. Define the following concomitants
\[ \begin{aligned} 
\Phi(2,4,2) & = \alpha_x^2 \, \beta_x^2 \, \brac{\alpha}{\beta}{u}^2 \\ 
\Phi(2,0,4) & = \brac{\alpha}{\beta}{u}^4 \\ 
\Phi(3,6,0) & = \alpha_x^2 \, \beta_x^2 \, \gamma_x^2 \, 
\brac{\alpha}{\beta}{\gamma}^2 \\
\Phi(3,3,3) & = \alpha_x \beta_x \gamma_x 
\brac{\alpha}{\beta}{\gamma} \, \brac{\alpha}{\beta}{u} \, 
\brac{\alpha}{\gamma}{u} \, \brac{\beta}{\gamma}{u} \\ 
\Phi(3,2,2) & = \gamma_x^2 \, \brac{\alpha}{\beta}{\gamma}^2 \, 
\brac{\alpha}{\beta}{u}^2 \\
\Phi(3,0,0) & = \brac{\alpha}{\beta}{\gamma}^4  \\ 
\Phi(3,0,6) & = \brac{\alpha}{\beta}{u}^2 \brac{\alpha}{\gamma}u^2 
\brac{\beta}{\gamma}{u}^2 \\ 
\Phi(4,0,2) & = \brac{\alpha}{\gamma}{\delta}^2 \brac{\beta}{\gamma}{\delta}^2 
\brac{\alpha}{\beta}{u}^2  \\
\Phi(4,1,3) & = \alpha_x \brac{\alpha}{\gamma}{\delta}^2
\brac{\beta}{\gamma}{u}^2 \brac{\alpha}{\beta}{\delta} 
\brac{\beta}{\delta}{u} \\ 
\Phi(4,4,0) & = \alpha_x \beta_x \gamma_x \delta_x \, 
\brac{\alpha}{\beta}{\gamma} \brac{\alpha}{\beta}{\delta} 
\brac{\beta}{\gamma}{\delta} \brac{\alpha}{\gamma}{\delta} \\ 
\Phi(4,2,4) & = \alpha_x \beta_x \, \brac{\alpha}{\gamma}{\delta} 
\brac{\beta}{\gamma}{\delta} \brac{\alpha}{\beta}{u}^2 
\brac{\gamma}{\delta}{u}^2 \\ 
\Phi_{I}(5,0,4) & = \brac{\alpha}{\beta}{\gamma}^4
\brac{\delta}{\epsilon}{u}^4 \\ 
\Phi_{II}(5,0,4) & = \brac{\alpha}{\beta}{\gamma}^2 
\brac{\delta}{\epsilon}{u}^2 \brac{\alpha}{\delta}{\epsilon}^2 
\brac{\beta}{\gamma}{u}^2 \\
\Phi(5,2,0) & = \alpha_x \beta_x \brac{\alpha}{\beta}{\gamma}^2 
\brac{\alpha}{\delta}{\epsilon} \brac{\beta}{\delta}{\epsilon} 
\brac{\gamma}{\delta}{\epsilon}^2 \\ 
\Phi(6,0,0) & = \brac{\alpha}{\beta}{\gamma}^2 
\brac{\delta}{\epsilon}{\zeta}^2 \brac{\alpha}{\epsilon}{\zeta}^2 
\brac{\beta}{\gamma}{\delta}^2 . 
\end{aligned} \] 

Now form the lists 
\[ \begin{aligned} 
\U{1} & = \{ \Phi(2,4,2), \Phi(2,0,4) \} \\ 
\U{2} & = \{ \Phi(3,6,0), \Phi(3,0,6), \Phi(3,3,3), 
            \Phi(3,2,2), \Phi(3,0,0)\} \\
\U{3} & = \{ \Phi(4,4,0), \Phi(4,2,4), \Phi(4,1,3), \Phi(4,0,2) \} \\
\U{4} & = \{ \Phi_{I}(5,0,4) - 3 \, \Phi_{II}(5,0,4), 
            \Phi(5,2,0) \} \\
\U{5} & = \{ 3 \, \Phi(6,0,0) - \Phi(3,0,0)^2 \}.
\end{aligned} \] 

If $\U{}$ is such a list, then $\U{}|_F = 0$ 
(resp. $\U{}|_F \neq 0$) means that all elements 
of ${\U{}}$ vanish at $F$ (resp.~at least one element 
is nonzero at $F$). 

With notation as above, our main theorem is the following: 
\begin{Theorem} \label{main.theorem}
For a ternary quartic $F$, 
\[ F \in W_s \iff \U{s}|_F = 0. \]  
\end{Theorem} 
In the sequel, it is frequently necessary to calculate plethysms 
and tensor products of $SL(V)$-representations, this was done 
using John Stembridge's SF package for Maple. All commutative 
algebra computations were done in Macaulay-2. 
\subsection{Case $s=1$} 
The locus $W_1$ is the quartic Veronese imbedding 
of $\P \, V^*$. It is well-known that its ideal is generated in 
degree $2$, and we have an exact sequence 
\[ \begin{array}{rccc}
0 & \ra H^0(\P^{14}, \II_{W_1}(2)) & \ra H^0(\P^{14}, \O_{\P^{14}}(2)) & \ra 
H^0(W_1, \O_{W_1}(2)) \\ 
  &                             &    ||        &  || \\ 
  &                             & S_2 (S_4 \, V)  &  S_8 \, V \end{array}
\] 
Decomposing $S_2(S_4 \, V)$, we see that 
$H^0(\P^{14}, \II_{W_1}(2))$ must be isomorphic to 
$ S_{(6,2)} V \oplus S_{(4,4)} V$. 
Henceforth we write the latter as $\{ 62, 44 \}$. 

To specify the inclusion $\{62\} \subseteq S_2 (S_4 \, V)$ is to specify 
a concomitant $\Phi(2,4,2)$. There is only one copy of $\{62\}$ inside 
$S_2 (S_4)$, hence there is a unique such $\Phi$ upto a constant. 
Now observe that 
$\alpha_x^2 \, \beta_x^2 \, \brac{\alpha}{\beta}{u}^2$ is a (legal) 
symbolic expression of the right degree, moreover it is not 
identically zero. 
(This is tantamount to checking that it is a nonzero element 
in the `bracket algebra' (see \cite{Sturmfels}). This was 
done in Macaulay-2.) Thus we have found $\Phi(2,4,2)$. The other 
concomitant $\Phi(2,0,4)$ is found in the same way, and this finishes 
the calculation for $s=1$. 

\begin{Remark} \rm 
In general, given $d,m,n$, it is possible to get all symbolic expressions 
of that degree by solving a system of Diophantine equations. However, 
in practice it is much easier to concoct such expressions by hand, 
especially if the multiplicity of $S_{m+n,n}$ in $S_d(S_4)$ is small. 
\end{Remark} 

\subsection{Case $s=2$.} \label{case.s2} 
Firstly we calculate the 
ideal $I_{W_2}$ by explicit elimination. Let 
\[ F = \sum\limits_{|I| = 4} a_I \, x^I, \quad 
L_i = p_{i0} \, x_0 + p_{i1} \, x_1 + p_{i2} \, x_2 \; \; \text{for $i=1,2$;}
\] 
where $a_I, p_{ij}$ are indeterminates. Write 
$F = L_1^4 + L_2^4$, and equate coefficients. We obtain polynomial 
expressions $a_I = f_I(p_{10}, \dots ,p_{22})$, defining a morphism 
\[ f: \complex \, [\{a_I\}] \lra \complex \, [\{p_{ij}\}] \] 
Then $I_{W_2}$ equals $\ker f $. The actual 
Macaulay-2 computation shows that all the ideal generators are 
in degree $3$, and $\dim \, (I_{W_2})_3 = 148$. 

The inclusion $W_2 \subseteq Y(2, \alpha_3) =Y$ implies 
$I_Y \subseteq I_{W_2}$. 
Now $Y$ is a rank variety in the sense of Porras \cite{Porras}, 
in particular it is reduced. 
(It is the locus of those $F$ which can be written as forms 
in only two variables by a change of coordinates.) The dimension 
of $Y$ is $6$ and its ideal sheaf is the image of the morphism 
(see [loc.cit.])
\[ \wedge^3 \alpha_3 : \wedge^3 (S_3 V) \otimes \O_{\P^{14}}(-3) 
\lra \O_{\P^{14}}. 
\] 
Decomposing $\wedge^3 (S_3)$, we deduce that 
$(I_Y)_3$ is the $120$-dimensional representation 
$\{63, 60, 42, 00\}$ inside 
\[ S_3 (S_4 V) = \{(12,0), (10,2),93,84,66,63,60,42,00\}. \] 
The quotient of the inclusion $(I_{Y})_3 \subseteq (I_{W_2})_3$ is 
a $28$-dimensional representation, and can only be $\{66\}$. 
Hence $(I_{W_2})_3 = \{66, 63, 60, 42, 00\}$. The concomitants 
are obtained as in the previous case. 

\begin{Remark} \rm 
The Gordan-N{\"o}ther theorem (see \cite[p.~234]{Olver}) 
implies that $F \in Y(2, \alpha_3)$ iff 
the Hessian of $F$ (which is $\Phi(3,6,0)$) identically vanishes. I do not 
know if the coefficients of $\Phi(3,6,0)$ define $Y$ scheme-theoretically. 
\end{Remark} 

\subsection{Case $s=3$.} Matters are greatly simplified 
because of the following lemma. 
\begin{Lemma} As schemes, $X_3 = Y(3, \alpha_2)$. 
\end{Lemma}
\demo 
As a first step, we show  that $X_3 = Y_\red(3, \alpha_2)$. 
Let $F \in Y_\red$. If $\rank \alpha_{1,F} \le 2$, then 
$F$ is a binary quartic in disguise, and then it has infinitely many 
apolar schemes of length $3$ (see \cite[\S 1.3]{IK}). 
If $\rank \alpha_{1,F} = 3$, then the existence of a 
length $3$ apolar scheme follows from \cite[Theorem 5.31]{IK}. 
(To summarise the situation, the Buchsbaum-Eisenbud structure theorem 
implies that $\ker \alpha_F$ is generated as an ideal by $3$ conics and 
$2$ quartics. The subideal generated by the $3$ conics
defines the apolar scheme.) In either case, $F \in X_3$. 
This shows that $X_3 = Y_\red(3, \alpha_2)$. 

Now $Y = Y(3, \alpha_2)$ is irreducible of dimension $8$, 
so equality holds in (\ref{eq.codim}). Hence $Y$ is Cohen-Macaulay, and 
has no embedded components. By the determintal formula, 
$\deg Y = 112$ which is the same as $\deg W_3$. Hence 
$Y$ must be reduced. \qed 

It follows from (\ref{res.jpw}) that $I_{W_3} = I_Y$ is generated 
in degree $4$, and 
\[ (I_{W_3})_4 = S_{2222}(S_2 \, V) = 
\{ 64, 43, 40,22\}. \] 
The concomitants are calculated as before. 

\subsection{Case $s=4$} This is similar to the previous case. 
\begin{Lemma} \sl 
As schemes, $X_4 = Y(4, \alpha_2)$. 
\end{Lemma} 
\demo 
Assume $F \in Y_\red(4, \alpha_2)$. If either 
$\rank \alpha_{1,F} \le 2$ or $\rank \alpha_{2,F} \le 3$, then 
$F \in X_3$ by the previous argument. Hence we may assume 
$\rank \alpha_{1,F} =3, \rank \alpha_{2,F} =4$. Then 
the two independent conics in $\ker \alpha_{2,F}$ define a 
complete intersection length $4$ scheme apolar to $F$. 
This shows that $X_4 = Y_\red(4, \alpha_2)$. The rest of 
the proof is similar to the previous lemma. 
\qed 

It follows that $I_{W_4}$ is generated in degree $5$ by 
$S_{22222}(S_2 \, V) = \{44, 20\}$. Now there are two copies 
of $\{44\}$ inside $S_5(S_4)$, hence a $2$-dimensional space 
of concomitants of degree $(5,0,4)$. A basis for this space is 
given by $\Phi_I(5,0,4), \Phi_{II}(5,0,4)$. Choose a typical 
form in $X_4$, say $F = x_0^4 + x_1^4 + x_2^4 + (x_0 + x_1 + x_2)^4$ 
and evaluate both concomitants at $F$. It is found that 
$\Phi_I - 3 \Phi_{II}$ vanishes identically on $F$. 

Similarly there are two copies of $\{20\}$ in $S_5(S_4)$. However it 
turns out that $\Phi(5,2,0)$ itself vanishes on $F$, so no 
linear combination is needed. 

\subsection{Case $s=5$.} The sextic invariant defining $W_5$ is 
called the catalecticant of ternary quartics. 
Decomposing $S_6(S_4 \, V)$, we see that it contains 
two copies of $\{00\}$. Now $\Phi(3,0,0)^2$ and $\Phi(6,0,0)$ 
generate this subspace, hence the  
catalecticant must be their linear combination. Specializing 
them at 
\[ 
F = x_0^4 + x_1^4 + x_2^4 + (x_0 + x_1 + x_2)^4 + (x_0 - x_1 + x_2)^4 , 
\] 
it is seen that $3 \, \Phi(6,0,0) - \Phi(3,0,0)^2$ identically 
vanishes at $F$. This completes the discussion of 
Theorem \ref{main.theorem}. \qed 

\subsection{A description of $W_2^\circ$} 
\label{w2open.subsection} 
In general $W_s^\circ$ is only expressible as a complicated boolean 
expression in closed sets, and it is not easy to characterize it 
algebraically. Here we attempt such a characterization for 
$s=2$. 

Let $F \in W_2 \setminus W_2^\circ$, then $F$ is apolar to  
a nonreduced length two subscheme $Z$ of $\P V^*$. Upto a change of 
coordinates, $I_Z = (y_0, y_1^2)$. This forces 
$F = c_1 \, x_2^4 + c_2 \, x_1 x_2^3$, for some constants $c_i$. 
Since $F$ has no apolar scheme of length one, $c_2 \neq 0$; so 
$F = x_2^3(\frac{c_1}{c_2} x_2 + x_1)$. 
Hence 
\[ W_2 \setminus W_2^\circ = 
\{ L_1^3L_2 : L_i \; \text{are linearly independent} \}. 
\] 
Now let 
\[ B = (W_2 \setminus W_2^\circ) \cup W_1 = 
\{ L_1^3L_2 : L_i \in V^* \}, 
\]  
which is an irreducible projective variety of dimension $4$. 
Geometrically, $B$ is the union of tangent planes to 
$W_1$. The inclusions $W_1 \subseteq B \subseteq W_2$ imply 
$I_{W_2} \subseteq I_B \subseteq I_{W_1}$. 

As in \S \ref{case.s2}, we calculate the generators 
of $I_B$ by explicit elimination. Its minimal resolution begins as 
\[ \begin{aligned} 
R(-3) \otimes M_8 \, \oplus R(-4) & \otimes M_{570} \, \oplus 
R(-5) \otimes M_{66} \ra \\ 
R(-2) & \otimes M_{15} \, \oplus R(-3) \otimes M_{56} \ra 
R \ra R/I_B \ra 0; 
\end{aligned} \] 
where $M_i$ is an $i$-dimensional $SL(V)$-representation. 
We need to identify $M_{15}$ and $M_{56}$. 
Since $(I_B)_2 \subseteq (I_{W_1})_2$, on dimensional grounds 
$M_{15} = \{44\}$. We have inclusions 
$(I_{W_2})_3 \subseteq (I_B)_3 \subseteq (I_{W_1})_3$, and the end terms 
are known. This forces 
\[ 
(I_B)_3 = \{84, 66, 63, 60, 42, 00\}. 
\] 
Now 
\[ M_8 \subseteq M_{15} \otimes R_1 = 
\{84, 63, 42, 21, 00 \}, \] 
hence $M_8 = \{21\}$. This implies that the submodule $\{84,63,42,00\}$ of 
$(I_B)_3$ is generated by $M_{15}$, so $M_{56}$ 
(the new generators in degree $3$) must be $\{66, 60 \}$. Define 
\[ {\mathcal V} = \{\Phi(2,0,4), \Phi(3,0,6), \Phi(3,6,0)\}, \] 
following the generators of $I_B$. 
Since $W_2^\circ = (W_2 \setminus B) \cup W_1$, we deduce that 
\begin{Proposition} \sl 
For a ternary quartic $F$, 
\[ F \in W_2^\circ \iff 
(\U{2}|_F = 0 \, \wedge \, {\mathcal V}|_F \neq 0) 
\vee (\U{1}|_F = 0). 
\] 
\end{Proposition} 
The general case does not seem so accessible, partly because 
there are a great many possibilities for the structure of a 
nonreduced scheme of length $s$. 

\bigskip 

{\small 
Acknowledgements: 
I am grateful to D.~Grayson and M.~Stillman (authors of Macaulay-2), 
as well as J.~Stembridge (author of the Symmetric Functions 
package for Maple). For information about obtaining these programs, see 
the respective websites: \\ `www.math.uiuc.edu/Macaulay2' and 
`www.math.umich.lsa.edu/$\sim$jrs'.}

\bibliographystyle{plain}
\bibliography{../../BIBLIO/ref1,../../BIBLIO/ref2,../../BIBLIO/ref3}

\vspace{1cm} 

\parbox{7cm}{ \small
Jaydeep V. Chipalkatti \\ 
Department of Mathematics \\ 
University of British Columbia \\ 
Vancouver B.C. V6T 1Z2,
Canada. \\
email: {\tt jaydeep@math.ubc.ca}}
\end{document}